\documentclass[11pt,a4paper,twoside]{article}
\usepackage{amssymb}
\usepackage{amsmath}
\usepackage{fancyhdr}
\usepackage{amsthm}
\usepackage{rotate}
\usepackage{graphicx}
\usepackage{float}
\usepackage[T1]{fontenc}

\usepackage{afterpage}  %

\newtheorem{theorem}{Theorem}

\newtheorem{example}{Example}

\theoremstyle{remark}

\begin{document}

\afterpage{\rhead[]{\thepage} \chead[\small I. I. Deriyenko
]{\small Configurations of conjugate permutations}
\lhead[\thepage]{} }

\begin{center}
\vspace*{2pt}
{\Large \textbf{Configurations of conjugate permutations}}\\[30pt]
{\large \textsf{\emph{Ivan I. Deriyenko}}}
\\[30pt]
{\sc Devoted to the memory of Valentin D. Belousov (1925-1988)}
\end{center}
{\bf Abstract.} {\footnotesize We describe some configurations of
conjugate permutations which may be used as a mathematical model
of some genetical processes and crystal growth. }
 \footnote{

\hspace*{2.5mm}flock, telomere, configuration.}

\section*{\centerline{1. Introduction}}\setcounter{section}{1}

Let $Q=\{1,2,3,\ldots,n\}$ be a finite set. The set of all
permutations of $Q$ will be denoted by $\mathbb{S}_n$. The
multiplication (composition) of permutations $\varphi$ and $\psi$
of $Q$ is defined as $\varphi\psi(x)=\varphi(\psi(x))$.
Permutations will be written in the form of cycles and cycles will
be separated by points, e.g.

$$
\varphi=\left(
\begin{array}{cccccc}
1 & 2 & 3 & 4 & 5 & 6\\
3 & 1 & 2 & 5 & 4 & 6
\end{array}
\right)=(123.45.6.)
$$

By a {\it type} of a permutation $\varphi\in\mathbb{S}_n$ we mean
the sequence
$$
C(\varphi)=\{l_1,l_2,\ldots,l_n\},
$$
where $l_i$ denotes the number of cycles of the length $i$.
Obviously,
$$
\sum_{i=1}^{n}i\cdot l_i=n\,.
$$

For example, for $\varphi=(132.45.6.)$ we have
$C(\varphi)=\{1,1,1,0,0,0\}$; for $\psi=(123456.)$ we obtain
$C(\psi)=\{0,0,0,0,0,1\}$.

As is well-known, two permutations $\varphi,\psi\in\mathbb{S}_n$
are {\it conjugate} if there exists a permutation
$\rho\in\mathbb{S}_n$ such that
\begin{equation}\label{e1}
\rho\varphi\rho^{-1}=\psi .
\end{equation}

\begin{theorem}\label{T1}{\rm (Theorem 5.1.3 in \cite{1})}
Two permutations are conjugated if and only if they have the same
type.\hfill$\Box{}$
\end{theorem}

In this short note we find all solutions of \eqref{e1}, i.e., for
a given $\varphi$ and $\psi$ we find all permutations $\rho$
satisfying this equation, and describe some graphs connected with
these solutions.

\section*{\centerline{2. Solutions of the equation (1)}}

Let's consider the equation \eqref{e1}. If
$\varphi=\psi=\varepsilon$, then as $\rho$ we can take any
permutation from $\mathbb{S}_n$. So, in this case \eqref{e1} has
$n!$ solutions.

If permutations $\varphi$ and $\psi$ are cyclic, then without loss
of generality, we can assume that
$$
\begin{array}{l}
\varphi=(1\,\varphi(1)\,\varphi^2(1)\,\varphi^3(1)\,\ldots\,\varphi^{n-1}(1).),\\[4pt]
\psi=(1\,\psi(1)\,\psi^2(1)\,\psi^3(1)\,\ldots\,\psi^{n-1}(1).),
\end{array}
$$
where $\varphi^0(1)=\varphi^n(1)=1$ and $\psi^0(1)=\psi^n(1)=1$.
In this case for $\rho_0$ defined by
\begin{equation}\label{e2}
\rho_0(\varphi^i(1))=\psi^i(1)=x_i, \ \ \ i=0,1,\ldots,n-1,
\end{equation}
we have
\[
\rho_0\varphi\rho_0^{-1}(x_i)=\rho_0\varphi\rho_0^{-1}(\psi^i(1))=\rho_0\varphi^{i+1}(1)=\psi^{i+1}(1)=\psi(\psi^i(1))=\psi(x_i),
\]
which shows that $\rho_0$ satisfies \eqref{e1}. Moreover, as is
not difficult to see, each permutation of the form
\begin{equation}\label{e3}
\rho=\rho_0\varphi^i, \ \ \ i=0,1,\ldots,n-1
\end{equation}
also satisfies this equation. There are no other solutions. So, in
this case we have $n$ different solutions.

In the general case when $\varphi$ and $\psi$ are decomposed into
cycles of the length $k_1,k_2,\ldots,k_r$, i.e.,
$$
\begin{array}{l}
\varphi=(a_{11}\thinspace a_{12}\ldots a_{1k_1})\ldots(a_{r1}\ldots a_{rk_r}), \\[3pt]
\psi=(b_{11}\thinspace b_{12}\ldots b_{1k_1})\ldots(b_{r1}\ldots
b_{rk_r}),
\end{array}
$$
the solution $\rho$, according to \cite{1}, has the form
\begin{equation}
\label{2} \beta=\left(
\begin{array}{cccccccc}
a_{11} & a_{12} & \ldots & a_{1k_1} & \ldots & a_{r1} & \ldots & a_{rk_r} \\
b_{11} & b_{12} & \ldots & b_{1k_1} & \ldots & b_{r1} & \ldots &
b_{rk_r}
\end{array}
\right),
\end{equation}
where the first row contains all elements of $\varphi$, the second
-- elements of $\psi$ written in the same order as in
decompositions of $\varphi$ and $\psi$ into cycles. Replacing in
$\varphi$ the cycle $(a_{11}\thinspace a_{12}\ldots a_{1k_1})$ by
$(a_{12}\thinspace a_{13}\ldots a_{1k_1}\thinspace a_{11})$ we
save the permutation $\varphi$ but we obtain a new $\rho$. Similar
to arbitrary cycles of $\varphi$ and $\psi$. In this way we obtain
all $\rho$ satisfying \eqref{e1}.

Let's observe that the cycle $(a_{11}\thinspace a_{12}\ldots
a_{1k_1})$ gives $k_1$ possibilities for the construction $\rho$.
From $m$ cycles of the length $k$ we can construct $m!\,k^m$
various $\rho$. So, in the case
$C(\varphi)=C(\psi)=\{l_1,l_2,\ldots,l_n\}$ we can construct
$$
N_{\varphi}=l_1!\cdot l_2!\cdot 2^{l_2}\cdot l_3!\cdot
3^{l_3}\cdot\ldots\cdot l_n!\cdot n^{l_n}
$$
various $\rho$.

\section*{\centerline{3. Configurations of conjugate permutations}}  

As is well-known, any permutation $\varphi$ of the set $Q$ of
order $n$ can be decomposed into $r\leqslant n$ cycles of the
length $k_1,k_2,\ldots,k_r$ with $k_1 + k_2 + \ldots + k_r = n$.
We denote this fact by
$$
Z=Z(\varphi)=[k_1,k_2,\ldots,k_r]
$$
and assume that $k_1\leqslant k_2\leqslant\ldots\leqslant k_r$.
$Z(\varphi)$ is called the {\it cyclic type} of $\varphi$. The set
of all permutations of the set $Q$ with the same cyclic type $Z_i$
is denoted by $F_i$ and is called a {\it flock}. Permutations
belonging to the same flock are conjugate (Theorem \ref{T1}). The
number of flocks $F_i\subset\mathbb{S}_n$ is equal to the number
of possible decompositions of $n$ into a sum of natural numbers.

In each flock we select one permutation $\sigma$ and call it a
{\it stem-permutation}. For simplicity we can assume that elements
of this permutation are written in the natural order.

\begin{example}\rm\label{Ex1}
Let's consider the set $Q=\{1,2,3,4,5\}$. The number $5$ has seven
decompositions into a sum of natural numbers, so the set of all
permutations of $Q$ has seven flocks. Below we present these
flocks and their stem-permutations.

\hspace*{16mm}$
\begin{array}{lll}
Z_1:5=5         &&\sigma=(12345.)\\
Z_2:5=1+4       &&\sigma=(1.2345.)\\
Z_3:5=2+3       &&\sigma=(12.345.)\\
Z_4:5=1+2+2     &&\sigma=(1.23.45.)\\
Z_5:5=1+1+3     &&\sigma=(1.2.345.)\\
Z_6:5=1+1+1+2   &&\sigma=(1.2.3.45.)\\
Z_7:5=1+1+1+1=1 &&\sigma=(1.2.3.4.5.)=\varepsilon
.\hspace{20mm}\Box{}
\end{array}$
\end{example}

Let's consider an arbitrary flock $F_i\subset\mathbb{S}_n$ and its
stem-permutation $\sigma$. For an arbitrary permutation
$\varphi_0\in F_i$ we define the sequence of permutations
$\varphi_0, \varphi_1,\varphi_2,\ldots$ by putting
\begin{equation}\label{e5}
\varphi_{k+1}=\varphi_k\sigma\varphi_k^{-1}.
\end{equation}
Obviously all $\varphi_k$ are in $F_i$. The set $F_i$ is finite,
so $\varphi_p=\varphi_s$ for some $p$ and $s$.

\unitlength1mm \hspace*{30mm}\begin{picture}(60,12)
 \put(-1,0){\circle{1}}\put(-2,-4){$\varphi_0$}
  \put(0,0){\vector(1,0){6}}
  \put(7,0){\circle{1}}\put(5,-4){$\varphi_1$}
   \put(8,0){\vector(1,0){6}}
   \put(15,0){\circle{1}}\put(13,-4){$\varphi_2$}
    \put(16,0){\vector(1,0){6}}
     \put(24,0){\circle*{0.6}}\put(27,0){\circle*{0.6}}\put(30,0){\circle*{0.6}}
    \put(32,0){\vector(1,0){6}}
    \put(39,0){\circle{1}}\put(37,-4){$\varphi_p$}\put(37,3){$\varphi_s$}
       \put(39.5,1){\vector(1,1){5}}\put(45,6.2){\circle{1}}
       \put(44.5,-6){\vector(-1,1){5}}\put(45,-6.2){\circle{1}}
       \put(46,6.5){\vector(1,0){6}}\put(53,6.2){\circle{1}}
       \put( 52,-6.5){\vector(-1,0){6}}\put(53,-6.2){\circle{1}}
    \put(54,6){\vector(3,-2){4}}
    \put(58,-3.5){\vector(-3,-2){4}}
    \put(58.5,-2){\circle*{0.6}}\put(58.5,0){\circle*{0.6}}\put(58.5,2){\circle*{0.6}}
\end{picture}

\medskip\medskip\medskip\medskip

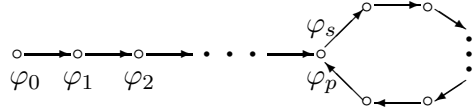
\centerline{\footnotesize{Fig. 1. The graph connected with the
sequence \eqref{e5}.}}

\medskip

The sequence $\varphi_1,\varphi_2,\varphi_3,\ldots$ can be
initiated by various $\varphi_0$ because for fixed $\varphi_1$ and
$\sigma$  the equation $\varphi_1=\varphi\sigma\varphi^{-1}$ has
many solutions.

Let's denote by $\Phi_k$ the set of all possible solutions of the
equation \eqref{e5}, where $\varphi_{k+1}$ and $\sigma$ are fixed.
Let
$$
\overline{\Phi}_k=\{\varphi\in\Phi_k\,:\,Z(\varphi)=Z(\sigma)\}.
$$
In the case when $\overline{\Phi}_k$ has only one element the
permutation $\varphi_{k+1}$ is called {\it simple}. If
$\overline{\Phi}_k$ is the empty set, then $\varphi_{k+1}$ is
called a {\it telomere} and is denoted by $\hat{\varphi}_{k+1}$.
In the corresponding oriented graph a telomere is a vertex which
is not preceded by another vertex.

\medskip

The following theorem is obvious.

\begin{theorem}
Let $\sigma$ be a stem-permutation of a flock $F_i$. If
$\varphi\in F_i$ is a telomere, then also
$\psi=\sigma\varphi\sigma^{-1}$ is a telomere. \hfill$\Box{}$
\end{theorem}

Two permutations $\varphi,\psi\in F_i\subset\mathbb{S}_n$ have the
same {\it configuration} $K$ if $\varphi_p=\psi_q$ for some
natural $p$ and $q$, where

$\begin{array}{lll}
& \ \ \ \ \ \ \ \ \ \ \ \ &\varphi_p=\varphi_{p-1}\sigma\varphi_{p-1}^{-1}\,,\,\ldots\,,\,\varphi_1=\varphi\sigma\varphi^{-1},\\[2mm]
&&\psi_q=\psi_{q-1}\sigma\psi_{q-1}^{-1}\,,\,\ldots\,,\,\psi_1=\psi\sigma\psi^{-1}
\end{array}
$

\noindent and $\sigma$ is a stem-permutation from $F_i$.

\section*{\centerline{4. A simple algorithm for determining
configurations}}

\noindent {\bf 1.} In a given flock $F_i$ we select a
stem-permutation $\sigma$ and one permutation
$\varphi_0\ne\sigma$. Using these two permutations and \eqref{e5}
we construct the sequence $\varphi_0, \varphi_1,\ldots,\varphi_l$,
where $\varphi_l\ne\varphi_s$ for all $0\leqslant s<l$ and
$\varphi_{l+1}=\varphi_t$ for some $0\leqslant t<l$. In this way
we obtain the graph\\

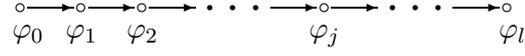
\unitlength1mm \hspace*{30mm}\begin{picture}(60,12)
 \put(-1,10){\circle{1}}\put(-2,6){$\varphi_0$}
  \put(0,10){\vector(1,0){6}}
  \put(7,10){\circle{1}}\put(5,6){$\varphi_1$}
   \put(8,10){\vector(1,0){6}}
   \put(15,10){\circle{1}}\put(13,6){$\varphi_2$}
    \put(16,10){\vector(1,0){6}}
     \put(24,10){\circle*{0.6}}\put(27,10){\circle*{0.6}}\put(30,10){\circle*{0.6}}
    \put(32,10){\vector(1,0){6}}
    \put(39,10){\circle{1}}\put(40,10){\vector(1,0){6}}\put(37,6){$\varphi_j$}
     \put(48,10){\circle*{0.6}}\put(51,10){\circle*{0.6}}\put(54,10){\circle*{0.6}}
    \put(56,10){\vector(1,0){6}}
    \put(63,10){\circle{1}}\put(62,6){$\varphi_l$}
\end{picture}

\noindent {\bf 2.} For each $\varphi_j$ from the above sequence,
from all solutions of the equation
$$
\rho\sigma\rho^{-1}=\varphi_j
$$
we select these solutions $\rho\ne\varphi_{j-1}$ which are in
$F_i$ and attach them to the previous solutions as immediately
preceding $\varphi_j$. In this way we obtain the configuration
$K=\{\varphi_0, \varphi_1,\ldots,\varphi_l,\rho_1,\rho_2,\ldots\}$
and the graph

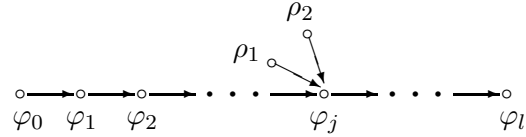
\unitlength1mm \hspace*{30mm}\begin{picture}(60,28)
 \put(-1,10){\circle{1}}\put(-2,6){$\varphi_0$}
  \put(0,10){\vector(1,0){6}}
  \put(7,10){\circle{1}}\put(5,6){$\varphi_1$}
   \put(8,10){\vector(1,0){6}}
   \put(15,10){\circle{1}}\put(13,6){$\varphi_2$}
    \put(16,10){\vector(1,0){6}}
     \put(24,10){\circle*{0.6}}\put(27,10){\circle*{0.6}}\put(30,10){\circle*{0.6}}
    \put(32,10){\vector(1,0){6}}
    \put(39,10){\circle{1}}\put(40,10){\vector(1,0){6}}\put(37,6){$\varphi_j$}
     \put(48,10){\circle*{0.6}}\put(51,10){\circle*{0.6}}\put(54,10){\circle*{0.6}}
    \put(56,10){\vector(1,0){6}}
    \put(63,10){\circle{1}}\put(62,6){$\varphi_l$}

    \put(32.1,14.2){\circle{1}}\put(32.53,14){\vector(2,-1){6}}\put(27,15){$\rho_1$}
    \put(36.8,18){\circle{1}}\put(37,17.5){\vector(1,-3){2}}\put(34,20.5){$\rho_2$}
\end{picture}

Next, for all new $\rho_k$ attached to $K$ we solve the equation
$\rho\sigma\rho^{-1}=\rho_k$ and attach to $K$ these solutions
$\rho^{\,\prime}\ne\rho_k$ which are in $F_i$. For this new
$\rho^{\prime}$ we solve the equation
$\rho\sigma\rho^{-1}=\rho^{\prime}$ and so on. Since $F_i$ is
finite after some steps we obtain a telomere which completes this
procedure.

\section*{\centerline{5. Examples}}

Now we give some examples. We will consider the set
$Q=\{1,2,3,4,5,6\}$ and its permutations. For simplicity we
consider the flock $F_1$ containing all cyclic permutations of $Q$
and select $\sigma=(123456.)$ as a stem-permutation of $F_1$.

\begin{example}\rm
If we choose $\varphi_0=(125634.)$, then, according to \eqref{e5},
we obtain
\[\begin{array}{lll}
&&\varphi_1=\varphi_0\sigma\varphi_0^{-1}=(163254.), \\[2pt]
&&\varphi_2=\varphi_1\sigma\varphi_1^{-1}=(143625.), \\[2pt]
&&\varphi_3=\varphi_2\sigma\varphi_2^{-1}=(163254.)=\varphi_1\, .
\end{array}
\]
Thus, the first step of our algorithm gives the configuration
$K=\{\varphi_0,\varphi_1,\varphi_2\}$.

Now, for each $\varphi_i\in K$ we solve the equation
$\rho\sigma\rho^{-1}=\varphi_i$ and add to $K$ all solutions
belonging to $F_1$.

The equation $\rho\sigma\rho^{-1}=\varphi_0$ is satisfied by the
permutation $\rho_0=(1.2.35.46.)$. So, according to \eqref{e3},
other solutions of this equation have the form

\[
\begin{array}{ll}
\varphi_{01}=\rho_0\sigma \ =(1.2.35.46.)(123456.)\ = (125436.),\\[3pt]
\varphi_{02}=\rho_0\sigma^2=(1.2.35.46.)(135.246.)=(15.26.3.4.),\\[3pt]
\varphi_{03}=\rho_0\sigma^3=(1.2.35.46.)(14.25.36.)=(165234.),\\[3pt]
\varphi_{04}=\rho_0\sigma^4=(1.2.35.46.)(153.264.) =(13.24.5.6.),\\[3pt]
\varphi_{05}=\rho_0\sigma^5=(1.2.35.46.)(165432.) \ =(145632.).
\end{array}
\]

From these solutions only $\varphi_{01},\varphi_{03},\varphi_{05}$
are in $F_1$. We attach these solutions to $K$ as the immediately
preceding $\varphi_0$.

Next, we consider the equation $\rho\sigma\rho^{-1}=\varphi_1$.
This equation has only one solution belonging to $F_1$. Since this
solution coincides with $\rho$, we do not obtain permutations
which should be added to $K$.

The equation $\rho\sigma\rho^{-1}=\varphi_2$ has only one solution
$\rho=(145236.)\ne\varphi_1$ belonging to $F_1$. We denote it by
$\varphi_4$ and add to $K$ as the solution immediately preceding
$\varphi_2$. At this instant we have the configuration
(uncomplete)
$$
K=\{\varphi_0,\varphi_1,\varphi_2,\varphi_{01},\varphi_{03},\varphi_{05},\varphi_4\}
$$
and the graph

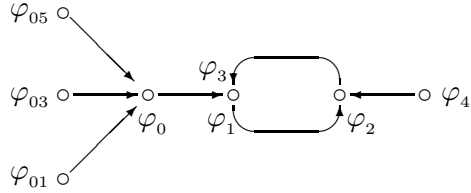
\unitlength1.4mm \hspace*{45mm}\begin{picture}(60,22)
 \put(-1,2){\circle{1}}\put(-6,2){$\varphi_{_{01}}$}
  \put(-0.3,2.7){\vector(1,1){6.2}}
  \put(-1,10){\circle{1}}\put(-6,9.5){$\varphi_{_{03}}$}
  \put(0,10){\vector(1,0){6}}
  \put(-1,17.8){\circle{1}}\put(-6,17.5){$\varphi_{_{05}}$}
  \put(-0.3,17.3){\vector(1,-1){6}}
  \put(7,10){\circle{1}}\put(6,7){$\varphi_{_0}$}
   \put(8,10){\vector(1,0){6}}
   \put(15,10){\circle{1}}\put(15,12){\vector(0,-1){1}}\put(12.5,7){$\varphi_{_1}$}
   \put(11.7,12){$\varphi_{_3}$}
    \put(32,10){\vector(-1,0){6}}
   \put(25,10){\circle{1}}\put(25.6,7){$\varphi_{_2}$}
   \put(20,11){\oval(10,5)[t]} \put(20,9){\oval(10,5)[b]}
   \put(25,8){\vector(0,1){1}}
   \put(33,10){\circle{1}}\put(34.5,9.5){$\varphi_{_4}$}
\end{picture}

\medskip

Further we will work with the permutations $\varphi_{01}$,
$\varphi_{03}$, $\varphi_{05}$, $\varphi_4$. Equations
$\rho\sigma\rho^{-1}=\varphi_{0i}$, $i=1,3,5$, do not have
solutions belonging to $F_i$. So, $\varphi_{01}$, $\varphi_{03}$,
$\varphi_{05}$ are telomeres. We denote them by
$\hat{\varphi}_{01}$, $\hat{\varphi}_{03}$, $\hat{\varphi}_{05}$.

The equation $\rho\sigma\rho^{-1}=\varphi_4$ has three solutions
belonging to $F_1$. Namely,
\[
\begin{array}{ll}
\varphi_{41}=\rho^{\prime}\sigma \ =(1.6.24.35.)(123456.)\ = (143256.),\\[3pt]
\varphi_{43}=\rho^{\prime}\sigma^3=(1.6.24.35.)(14.25.36.)=(123654.),\\[3pt]
\varphi_{45}=\rho^{\prime}\sigma^5=(1.6.24.35.)(165432.)\
=(163452.).
\end{array}
\]

Since equations $\rho\sigma\rho^{-1}=\varphi_{4j}$, $j=1,3,5$, do
not have solutions belonging to $F_1$, $\varphi_{41}$,
$\varphi_{43}$, $\varphi_{45}$ are telomeres.

Summarizing the above we obtain the configuration
$$
K=\{\varphi_0,\varphi_1,\varphi_2,\hat{\varphi}_{01},\hat{\varphi}_{03},
\hat{\varphi}_{05},\varphi_4,\hat{\varphi}_{41},\hat{\varphi}_{43},\hat{\varphi}_{45}\}
$$
and the graph

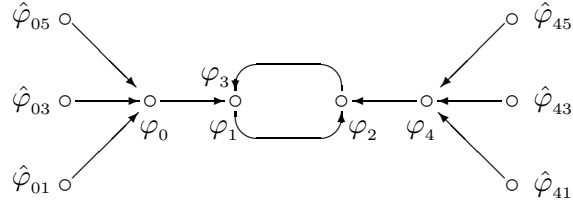
\unitlength1.4mm \hspace*{35mm}\begin{picture}(60,22)
 \put(-1,2){\circle{1}}\put(-6,2){$\hat{\varphi}_{_{01}}$}
  \put(-0.3,2.7){\vector(1,1){6.2}}
  \put(-1,10){\circle{1}}\put(-6,9.5){$\hat{\varphi}_{_{03}}$}
  \put(0,10){\vector(1,0){6}}
  \put(-1,17.8){\circle{1}}\put(-6,17.5){$\hat{\varphi}_{_{05}}$}
  \put(-0.3,17.3){\vector(1,-1){6}}
  \put(7,10){\circle{1}}\put(6,7){$\varphi_{_0}$}
   \put(8,10){\vector(1,0){6}}
   \put(15,10){\circle{1}}\put(15,12){\vector(0,-1){1}}\put(12.5,7){$\varphi_{_1}$}
   \put(11.7,12){$\varphi_{_3}$}
    \put(32,10){\vector(-1,0){6}}
   \put(25,10){\circle{1}}\put(25.6,7){$\varphi_{_2}$}
   \put(20,11){\oval(10,5)[t]} \put(20,9){\oval(10,5)[b]}
   \put(25,8){\vector(0,1){1}}
   \put(33,10){\circle{1}}\put(31,7){$\varphi_{_4}$}
   \put(41,2){\circle{1}}\put(43,1.5){$\hat{\varphi}_{_{41}}$}
  \put(40.3,2.7){\vector(-1,1){6.2}}
  \put(41,10){\circle{1}}\put(43,9.5){$\hat{\varphi}_{_{43}}$}
  \put(40,10){\vector(-1,0){6}}
  \put(41,17.8){\circle{1}}\put(43,17.5){$\hat{\varphi}_{_{45}}$}
  \put(40.3,17.3){\vector(-1,-1){6}}

\end{picture}
\end{example}

\begin{example}\rm Using the same flock $F_1$ and the same $\sigma$
but selecting another $\varphi_0$ we can obtain another
configuration. For example by selecting $\varphi_0=(162435.)$ we
obtain the configuration $K_2$ presented by the following graph:

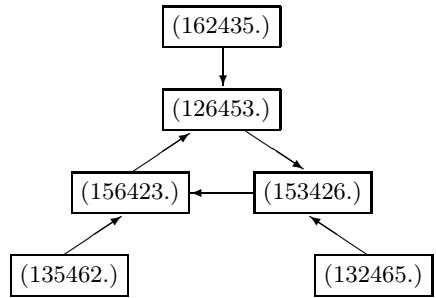
\unitlength1mm \hspace*{20mm}
\begin{picture}(60,44)
 \put(30,33){\fbox{\footnotesize{(162435.)}}}
 \put(38,31){\vector(0,-1){5}}
 \put(30,22){\fbox{\footnotesize{(126453.)}}}
 \put(42,11.82){\vector(-1,0){8.3}}
 \put(26,14.6){\vector(3,2){7.6}}\put(41,20){\vector(3,-2){7.6}}
 \put(18,11){\fbox{\footnotesize{(156423.)}}}\put(42,11){\fbox{\footnotesize{(153426.)}}}
 \put(17,3.6){\vector(3,2){7.6}}\put(57,3.6){\vector(-3,2){7.6}}
 \put(10,0){\fbox{\footnotesize{(135462.)}}}\put(50,0){\fbox{\footnotesize{(132465.)}}}
\end{picture}
\end{example}

\bigskip\medskip

\noindent {\bf Remark.} The flock $F_1$ has six configurations:

$\bullet$ \ $K_1$ and $K_2$ are described in the above examples,

$\bullet$ \ $K_3$ induced by $\varphi_0=(125643.)$ contains $18$
permutations,

$\bullet$ \ $K_4$ induced by $\varphi_0=(135624.)$ contains $42$
permutations,

$\bullet$ \ $K_5$ induced by $\varphi_0=(136245.)$ contains $42$
permutations,

$\bullet$ \ $K_6$ has only two permutations: $\sigma$ and
$\sigma^{-1}$.

\smallskip
\noindent Flocks $K_4$ and $K_5$ are isomorphic as graphs.

\smallskip
The set $\mathbb{S}_6$ is divided into $11$ flocks.

\medskip

The author does'nt know a general method that would allow to
determine the number of configurations in each flock. Neither does
he know how to quickly find a telomere using stem-permutations. It
is also unknown how to check if two telomeres belong to the same
configuration.

\section*{\centerline{6. Conclusions}}

The results shown were inspired by some research in genetics. Some
termino\-logy (stem-permutation, telomere) was also drawn from
genetics. The author thinks that the described method of
configuration can be effectively used in chemistry in researching
growth of crystals.

\noindent \footnotesize{\rightline{Received \ May 8, 2010}

\medskip\noindent
Department of Higher Mathematics and Informatics, Kremenchuk State
Polytechnic University,
20 Pervomayskaya str, 39600 Kremenchuk, Ukraine\\
E-mail: ivan.deriyenko@gmail.com}

\end{document}